\documentstyle[12pt,amssymb]{article}

\newcommand{\Label}{\label}

\begin{document}

\baselineskip15pt

\newtheorem{definition}{Definition $\!\!$}[section]
\newtheorem{prop}[definition]{Proposition $\!\!$}
\newtheorem{lem}[definition]{Lemma $\!\!$}
\newtheorem{corollary}[definition]{Corollary $\!\!$}
\newtheorem{theorem}[definition]{Theorem $\!\!$}
\newtheorem{example}[definition]{\it Example $\!\!$}
\newtheorem{remark}[definition]{Remark $\!\!$}

\newcommand{\nc}[2]{\newcommand{#1}{#2}}
\newcommand{\rnc}[2]{\renewcommand{#1}{#2}}

\nc{\Section}{\setcounter{definition}{0}\section}
\renewcommand{\theequation}{\thesection.\arabic{equation}}
\newcounter{c}
\renewcommand{\[}{\setcounter{c}{1}$$}
\newcommand{\etyk}[1]{\vspace{-7.4mm}$$\begin{equation}\Label{#1}
\addtocounter{c}{1}}
\renewcommand{\]}{\ifnum \value{c}=1 $$\else \end{equation}\fi}


\nc{\bpr}{\begin{prop}}
\nc{\bth}{\begin{theorem}}
\nc{\ble}{\begin{lem}}
\nc{\bco}{\begin{corollary}}
\nc{\bre}{\begin{remark}}
\nc{\bex}{\begin{example}}
\nc{\bde}{\begin{definition}}
\nc{\ede}{\end{definition}}
\nc{\epr}{\end{prop}}
\nc{\ethe}{\end{theorem}}
\nc{\ele}{\end{lem}}
\nc{\eco}{\end{corollary}}
\nc{\ere}{\hfill\mbox{$\Diamond$}\end{remark}}
\nc{\eex}{\end{example}}
\nc{\epf}{\hfill\mbox{$\Box$}}
\nc{\ot}{\otimes}
\nc{\bsb}{\begin{Sb}}
\nc{\esb}{\end{Sb}}
\nc{\ct}{\mbox{${\cal T}$}}
\nc{\ctb}{\mbox{${\cal T}\sb B$}}
\nc{\bcd}{\[\begin{CD}}
\nc{\ecd}{\end{CD}\]}
\nc{\ba}{\begin{array}}
\nc{\ea}{\end{array}}
\nc{\bea}{\begin{eqnarray}}
\nc{\eea}{\end{eqnarray}}
\nc{\be}{\begin{enumerate}}
\nc{\ee}{\end{enumerate}}
\nc{\beq}{\begin{equation}}
\nc{\eeq}{\end{equation}}
\nc{\bi}{\begin{itemize}}
\nc{\ei}{\end{itemize}}
\nc{\kr}{\mbox{Ker}}
\nc{\te}{\!\ot\!}
\nc{\pf}{\mbox{$P\!\sb F$}}
\nc{\pn}{\mbox{$P\!\sb\nu$}}
\nc{\bmlp}{\mbox{\boldmath$\left(\right.$}}
\nc{\bmrp}{\mbox{\boldmath$\left.\right)$}}
\rnc{\phi}{\mbox{$\varphi$}}
\nc{\LAblp}{\mbox{\LARGE\boldmath$($}}
\nc{\LAbrp}{\mbox{\LARGE\boldmath$)$}}
\nc{\Lblp}{\mbox{\Large\boldmath$($}}
\nc{\Lbrp}{\mbox{\Large\boldmath$)$}}
\nc{\lblp}{\mbox{\large\boldmath$($}}
\nc{\lbrp}{\mbox{\large\boldmath$)$}}
\nc{\blp}{\mbox{\boldmath$($}}
\nc{\brp}{\mbox{\boldmath$)$}}
\nc{\LAlp}{\mbox{\LARGE $($}}
\nc{\LArp}{\mbox{\LARGE $)$}}
\nc{\Llp}{\mbox{\Large $($}}
\nc{\Lrp}{\mbox{\Large $)$}}
\nc{\llp}{\mbox{\large $($}}
\nc{\lrp}{\mbox{\large $)$}}
\nc{\lbc}{\mbox{\Large\boldmath$,$}}
\nc{\lc}{\mbox{\Large$,$}}
\nc{\Lall}{\mbox{\Large$\forall\;$}}
\nc{\bc}{\mbox{\boldmath$,$}}
\rnc{\epsilon}{\varepsilon}
\rnc{\ker}{\mbox{\em Ker}}
\nc{\ra}{\rightarrow}
\nc{\ci}{\circ}
\nc{\cc}{\!\ci\!}
\nc{\T}{\mbox{\sf T}}
\nc{\can}{\mbox{\em\sf T}\!\sb R}
\nc{\cnl}{$\mbox{\sf T}\!\sb R$}
\nc{\lra}{\longrightarrow}
\nc{\M}{\mbox{Map}}
\rnc{\to}{\mapsto}
\nc{\imp}{\Rightarrow}
\rnc{\iff}{\Leftrightarrow}
\nc{\ob}{\mbox{$\Omega\sp{1}\! (\! B)$}}
\nc{\op}{\mbox{$\Omega\sp{1}\! (\! P)$}}
\nc{\oa}{\mbox{$\Omega\sp{1}\! (\! A)$}}
\nc{\inc}{\mbox{$\,\subseteq\;$}}
\nc{\de}{\mbox{$\Delta$}}
\nc{\spp}{\mbox{${\cal S}{\cal P}(P)$}}
\nc{\dr}{\mbox{$\Delta_{R}$}}
\nc{\dsr}{\mbox{$\Delta_{\cal R}$}}
\nc{\m}{\mbox{m}}
\nc{\0}{\sb{(0)}}
\nc{\1}{\sb{(1)}}
\nc{\2}{\sb{(2)}}
\nc{\3}{\sb{(3)}}
\nc{\4}{\sb{(4)}}
\nc{\5}{\sb{(5)}}
\nc{\6}{\sb{(6)}}
\nc{\7}{\sb{(7)}}
\nc{\hsp}{\hspace*}
\nc{\nin}{\mbox{$n\in\{ 0\}\!\cup\!{\Bbb N}$}}
\nc{\al}{\mbox{$\alpha$}}
\nc{\bet}{\mbox{$\beta$}}
\nc{\ha}{\mbox{$\alpha$}}
\nc{\hb}{\mbox{$\beta$}}
\nc{\hg}{\mbox{$\gamma$}}
\nc{\hd}{\mbox{$\delta$}}
\nc{\he}{\mbox{$\varepsilon$}}
\nc{\hz}{\mbox{$\zeta$}}
\nc{\hs}{\mbox{$\sigma$}}
\nc{\hk}{\mbox{$\kappa$}}
\nc{\hm}{\mbox{$\mu$}}
\nc{\hn}{\mbox{$\nu$}}
\nc{\la}{\mbox{$\lambda$}}
\nc{\hl}{\mbox{$\lambda$}}
\nc{\hG}{\mbox{$\Gamma$}}
\nc{\hD}{\mbox{$\Delta$}}
\nc{\th}{\mbox{$\theta$}}
\nc{\Th}{\mbox{$\Theta$}}
\nc{\ho}{\mbox{$\omega$}}
\nc{\hO}{\mbox{$\Omega$}}
\nc{\hp}{\mbox{$\pi$}}
\nc{\hP}{\mbox{$\Pi$}}
\nc{\bpf}{{\it Proof.~~}}
\nc{\slq}{\mbox{$A(SL\sb q(2))$}}
\nc{\fr}{\mbox{$Fr\llp A(SL(2,\IC))\lrp$}}
\nc{\slc}{\mbox{$A(SL(2,\IC))$}}
\nc{\af}{\mbox{$A(F)$}}
\rnc{\widetilde}{\tilde}
\nc{\suq}{\mbox{$A(SU_q(2))$}}
\nc{\asq}{\mbox{$A(S_q^2)$}}
\nc{\tasq}{\mbox{$\widetilde{A}(S_q^2)$}}

\def\esl{{\mbox{$E\sb{\frak s\frak l (2,{\Bbb C})}$}}}
\def\esu{{\mbox{$E\sb{\frak s\frak u(2)}$}}}
\def\zf{{\mbox{${\Bbb Z}\sb 4$}}}
\def\zt{{\mbox{$2{\Bbb Z}\sb 2$}}}
\def\ox{{\mbox{$\Omega\sp 1\sb{\frak M}X$}}}
\def\oxh{{\mbox{$\Omega\sp 1\sb{\frak M-hor}X$}}}
\def\oxs{{\mbox{$\Omega\sp 1\sb{\frak M-shor}X$}}}
\def\Fr{\mbox{Fr}}

\def\gal{-Galois extension}
\def\hge{Hopf-Galois extension}
\def\qg{quantum group}
\def\aqpb{algebraic quantum principal bundle}
\def\ta{\tilde a}
\def\tb{\tilde b}
\def\tc{\tilde c}
\def\td{\tilde d}
\def\st{\stackrel}
\def\<{\langle}
\def\>{\rangle}

\def\inbar{\,\vrule height1.5ex width.4pt depth0pt}
\def\IC{{\Bbb C}}
\def\IZ{{\Bbb Z}}
\def\IN{{\Bbb N}}


\title{
\vspace*{-15mm}{\Large\bf 
ALGEBRAIC \boldmath $K_0$ OF THE QUANTUM SPHERE AND 
NONCOMMUTATIVE INDEX THEOREM}
}
\author{
{\normalsize\sc Piotr M.~Hajac}\thanks{  
On leave from:
Department of Mathematical Methods in Physics, 
Warsaw University, ul.~Ho\.{z}a 74, Warsaw, \mbox{00--682~Poland}.
http://info.fuw.edu.pl/KMMF/ludzie\underline{~~}ang.html 
(E-mail: pmh@fuw.edu.pl)
}\\ 
\normalsize International School for Advanced Studies,
\vspace*{-0mm}\\
\normalsize Via Beirut 2--4, 34013 Trieste, Italy.\\
\vspace*{-0mm}
\normalsize http://www.damtp.cam.ac.uk/user/pmh33
}
\date{}
\maketitle

\vspace*{-5mm}
\begin{abstract}
The Noncommutative Index Theorem is used to prove that  the Chern character of 
quantum Hopf line bundles over the standard Podle\'s quantum
sphere equals the winding number of the representations defining these bundles. 
This result gives an estimate of the positive cone of the {\em algebraic} $K_0$ 
of the standard quantum sphere.
\end{abstract}

\section*{Introduction}

The goal of this paper is to compute the Chern character of  quantum Hopf
line bundles. We do it within the framework of the Hopf-Galois theory of
extensions of rings and Chern-Connes pairing between the 
cyclic cohomology and $K$-theory. We view the quantum two-sphere 
as the base space of the quantum principal Hopf fibration described algebraically 
as a \hge.
Noncommutative line bundles over the quantum sphere are constructed as left
and right projective bimodules over the coordinate ring of the quantum sphere.
They are the bimodules of colinear mappings indexed by the winding
number of one dimensional representations of~$U(1)$.  The Chern character of
the left projective module of such a quantum line bundle is proved to
coincide with the winding number of the representation defining the
bundle. The Chern character of the corresponding right
projective modules is shown to be equal to the opposite of its left counterpart,
whence it coincides with minus the winding number. 

In the following section, we recall a definition of a \hge\ and general construction
of the bimodule of an associated quantum vector bundle. We also recall 
the construction of $SL_q(2)$, quantum principal Hopf fibration, quantum Hopf
line bundles and the cyclic cohomology classes of the quantum sphere that pair 
non-trivially with $K$-theory.
 In Section~2 we extend the relevant considerations in~\cite{hm98}, and 
compute the pairing between the Chern cyclic cocycle
of the quantum sphere and the left projector matrices of quantum Hopf line bundles
for an arbitrary winding number.
This computation relies on the integrality of the pairing, which is implied by
the Noncommutative Index Theorem, and yields the winding number, as expected. 
Then we argue that the pairing with the right
projector matrices equals minus the pairing with the corresponding left matrices.
We conclude by noting that  the image of the positive cone of the algebraic
$K_0$ of the quantum sphere under the pairing with cyclic cohomology 
contains~$\IN\times\IZ$.

To focus attention and access straightforwardly the $C^*$-algebraic framework,
we work over~$\IC$.  We assume that $q$ is a non-zero element in $\IC$
that is not a root of~1. We tacitly assume the compact $*$-structure 
making $SL_q(2)$ into~$SU_q(2)$. However, we do not write it out explicitly, as it is 
an additional structure whose existence is crucial but a concrete form not
directly relevant to the considerations presented herein.
 For an introduction to quantum groups (comprehensive discussion of $SL_q(2)$ 
included)  we refer to \cite{k-ch95}.
For a concise and motivating treatment of the Hopf-Galois theory see~\cite{s-hj94}.
The Chern-Connes pairing and Noncommutative Index Theorem are elaborated
upon in~\cite{c-a94,l-jl97}, and \cite{c-a88} contains a compact account of
the matter. A brief note on classical Hopf line bundles within the general
context of $K$-theory and Noncommutative Geometry can be found 
in~\cite[p.101]{m-j95}. The generalisation to the non-standard Podle\'s quantum spheres
of the Chern-Connes pairing calculated in \cite{hm98} is carried out in~\cite{bm98}.

\section{Preliminaries}

In this section we recall basic definitions and known results used in the sequel.
We begin with a definition of a Hopf-Galois extension. \hge s  describe quantum
principal bundles the same way Hopf algebras describe \qg s. Here
a Hopf algebra $H$ plays the role of the algebra of functions on the
structure group, and the total space of 
a bundle is replaced by an  $H$-comodule algebra~$P$. 
\bde\Label{hgdef} 
Let $H$ be a Hopf algebra, $P$ be a right $H$-comodule algebra with multiplication
$m_P$ and coaction \dr, and 
$B:=P\sp{co H}:=\{p\in P\,|\; \dr\ p=p\ot 1\}$ the subalgebra of coinvariants. 
We say that P is an  {\em $H$-Galois extension} of $B$ iff the canonical 
left $P$-module right $H$-comodule map 
$
\chi:=(m\sb P\ot id)\circ (id\ot\sb B \dr )\, :\; P\ot\sb B P\lra P\ot H
$
is bijective. 
\ede

A natural next step is to consider
associated quantum vector bundles. More precisely, what we need here is a
replacement for the module of sections of an associated vector bundle. In the
classical case such sections can be equivalently described  as ``functions
of type $\varrho$" from the total space of a principal bundle to a vector 
space. We follow this construction in the quantum case by considering $B$-bimodules
of colinear maps (linear maps that preserve the comodule structure)
$\mbox{Hom}_\rho(V,P)$ associated with 
an $H$\gal\ $B\inc P$ via a corepresentation $\rho:V\ra V\ot H$ 
(see~\cite{d-m96,d-m97}). Under certain reasonable assumptions, these bimodules
are always left and right finitely generated projective. Thus we remain within
the paradigm of the Serre-Swan theorem.

Let us now exemplify the foregoing concepts. Recall first that $A(SL_q(2))$ is a 
Hopf algebra over $\IC$
generated by $1,\, a,\, b,\, c,\, d$ satisfying the 
following relations: 
\bea \Label{comm}
&& ab=q^{-1}ba~,~~ac=q^{-1}ca~,~~bd=q^{-1}db~,~~bc=cb~,~~cd=q^{-1}dc~, \nonumber \\ 
&& ad-da=(q^{-1}-q)bc~,~~ad-q^{-1}bc=da-qbc=1\, ,  
\eea
where $q\in\IC\setminus\{0\}$. 
The comultiplication $\Delta$, counit $\varepsilon$, and antipode
$S$ of \slq\ are defined by the following formulas:
\[ 
\Delta\pmatrix{a & b \cr c & d \cr}=
\pmatrix{a\ot1 & b\ot1 \cr c\ot1 & d\ot1 \cr} 
\pmatrix{1\ot a & 1\ot b \cr 1\ot c & 1\ot d \cr},
\]\[
\varepsilon\pmatrix{\alpha & \beta \cr \gamma & \delta\cr}=
\pmatrix{ 1 & 0 \cr 0 & 1 \cr},~~~
S\pmatrix{a & b \cr c & d \cr}=
\pmatrix{d &-qb \cr -q^{-1}c & a \cr}.
\]
From here we can proceed to the construction of the standard quantum sphere of 
Podle\'s and the quantum principal Hopf fibration.
The standard quantum sphere is singled out among the
principal series of Podle\'s quantum spheres by the property that it can be
constructed as a quantum quotient space~\cite{p-p87}. In algebraic terms it
means that its coordinate ring can be obtained as the subalgebra of coinvariants
of a comodule algebra. To carry out this construction, first we need the
right coaction on \slq\ of the commutative and cocommutative Hopf algebra 
$\IC[z,z^{-1}]$ generated by the grouplike element $z$ and its inverse. This
Hopf algebra can be obtained as the quotient of \slq\ by the Hopf ideal generated
by the off-diagonal generators $b$ and~$c$. Identifying the image of $a$ and
$d$ under the Hopf algebra surjection $\pi:\slq\ra\IC[z,z^{-1}]$ with $z$ and 
$z^{-1}$ respectively, we can describe the right coaction $\dr:=(id\ot\pi)\ci\hD$
by the formula:
\[
\dr\pmatrix{a & b \cr c & d \cr}=
\pmatrix{a\ot z & b\ot z^{-1} \cr c\ot z & d\ot z^{-1} \cr}.
\]
We call the subalgebra of coinvariants defined by this coaction the coordinate
ring of the (standard) {\em quantum sphere}, and denote it by~\asq. Using general
tools of the Hopf-Galois theory (e.g., \cite[Theorem~I]{s-hj90}), it is 
straightforward to prove that $\asq\inc\slq$ is a $\IC[z,z^{-1}]$\gal. 
 We refer to the quantum principal
bundle given by this \hge\ as the
{\em quantum principal Hopf fibration}. (An $SO_q(3)$ version of this quantum 
fibration was studied in~\cite{bm93}.)

Now we need to define quantum
Hopf line bundles associated to the just described Hopf $q$-fibration and
provide their projector matrices.
\bde\Label{line}
Let $\rho_\mu:\IC[z,z^{-1}]\ra \IC\ot \IC[z,z^{-1}]$, $\rho_\mu(1)=1\ot z^{-\mu}$, 
$\mu\in\IZ$, be a one-dimensional corepresentation of $\IC[z,z^{-1}]$.
We call the \asq-bimodule of colinear maps $\mbox{\em Hom}_{\rho_\mu}(\IC,\slq)$
the (bimodule of) {\em quantum Hopf line bundle of winding number~$\mu$}.
\ede
We deal here with one-dimensional corepresentations, so that we can identify colinear
maps with their value at~1. We have
\[
\mbox{Hom}_{\rho_\mu}(\IC,\slq)\tilde{=}\{p\in\slq\ |\;\dr p=p\ot z^{-\mu}\}=:P_\mu
\]
as \asq-bimodules. With the help of the PBW basis 
$a^kb^lc^m,\;b^pc^rd^s,\; k,l,m,p,r,s\in\IN_0,\; k>0$ of \slq,
one can show that 
\begin{eqnarray*}\Label{osum}
&&
P_\mu=\left\{
\begin{array}{ll}
\sum_{k=0}^{-\mu}\asq\ a^{-\mu-k}c^{k}=\sum_{k=0}^{-\mu}a^{-\mu-k}c^{k}\asq
& \mbox{for $\mu\leq 0$}\\
\sum_{k=0}^{\mu}\asq\ b^{k}d^{\mu-k}=\sum_{k=0}^{\mu}b^{k}d^{\mu-k}\asq
& \mbox{for $\mu\geq 0$}, 
\end{array}
\right.
\end{eqnarray*}
and $\slq=\bigoplus_{\mu\in\IZ}P_\mu$ (cf.~\cite[(1.10)]{mmnnu91}).
Since the goal of this paper is to compute the Chern-Connes pairing 
\cite[p.224]{c-a94} between
 quantum Hopf line bundles and the cyclic cohomology of the standard quantum sphere,
we need explicit formulas for the projector matrices of the former and generators
of the latter. To this end recall first that if $xuv=vu$, then 
$(u+v)^n=\sum_{k=0}^n{\scriptsize\pmatrix{n\cr k\cr}}_{\!x}u^kv^{n-k}$, where
\[
{\scriptsize\pmatrix{n\cr k\cr}}_{\!x}=
\frac{(x-1)...(x^n-1)}{(x-1)...(x^k-1)(x-1)...(x^{n-k}-1)}
\]
 are the {\em $x$-binomial coefficients} (e.g., see Section~IV.2 in~\cite{k-ch95}).
Now, following~\cite{hm98}, put 
\begin{eqnarray}\Label{ekl}
&&
(e_\mu)_{kl}=\left\{
\begin{array}{ll}
a^{-\mu-k}c^{k}{\scriptsize\pmatrix{-\mu\cr l\cr}}_{\!q^2}(-q)^lb^ld^{-\mu-l}
& \mbox{for $\mu\leq 0$}\\
b^{k}d^{\mu-k}{\scriptsize\pmatrix{\mu\cr l\cr}}_{\!q^2}(-q)^{-l}a^{\mu-l}c^l
& \mbox{for $\mu\geq 0$}. 
\end{array}
\right.
\end{eqnarray}
Then, for any $\mu\in\IZ$, $e_\mu\in M_{|\mu|+1}(\asq)$, $e_\mu^2=e_\mu$, and 
$\asq^{|\mu|+1}e_\mu$ is isomorphic to $P_\mu$ as a left \asq-module~\cite{hm98}. 
Similarly, put 
\begin{eqnarray}\Label{flk}
&&
(f_\mu)_{lk}=\left\{
\begin{array}{ll}
{\scriptsize\pmatrix{-\mu\cr l\cr}}_{\!q^2}(-q)^{-l}b^ld^{-\mu-l}a^{-\mu-k}c^{k}
& \mbox{for $\mu\leq 0$}\\
{\scriptsize\pmatrix{\mu\cr l\cr}}_{\!q^2}(-q)^{l}a^{\mu-l}c^lb^{k}d^{\mu-k}
& \mbox{for $\mu\geq 0$}. 
\end{array}
\right.
\end{eqnarray}
Then again, for any $\mu\in\IZ$, $f_\mu\in M_{|\mu|+1}(\asq)$, $f_\mu^2=f_\mu$, and 
$f_\mu\asq^{|\mu|+1}$ is isomorphic to $P_\mu$ as a right \asq-module~\cite{hm98}. 

To obtain the desired pairing, we need 
to evaluate  appropriate even cyclic cocycles on the left and right
projector matrices provided above.  The positive even cyclic cohomology
$HC^{2n}(\asq)$, $n>0$, is the image of the periodicity operator applied to
$HC^{0}(\asq)$. In degree zero it is given by two generators
(cohomologically non-trivial cyclic cocycles) and the kernel of the periodicity 
operator~\cite[p.174]{mnw91}.
Since the pairing is compatible with the action of 
the periodicity
operator, it is completely determined by the aforementioned two cyclic
0-cocycles (traces). 
(Everything else either pairs with $K_0(\asq)$ in the same
way or trivially.) These traces are explicitly provided 
in~\cite{mnw91}. One of them, denoted by $\tau^0$, is simply the restriction
to \asq\ of the counit map of~\slq. It can be argued that this trace detects
 the ``rank" of our quantum vector bundles~\cite[Remark~3.4]{hm98}.
The other trace is given by the following adaptation of \cite[(4.4)]{mnw91} to our 
special case of the
standard Podle\'s quantum sphere: 
\bea\Label{tr}
&&
\tau^1((ab)^m\hz^n)=\left\{
\begin{array}{ll}
(1-q^{2n})^{-1} & \mbox{for $n>0$, $m=0$,}\\
0 & \mbox{otherwise}, 
\end{array}
\right.
\nonumber\\ && \ \nonumber\\ &&
\tau^1((cd)^m\hz^n)=\left\{
\begin{array}{ll}
(1-q^{2n})^{-1} & \mbox{for $n>0$, $m=0$,}\\
0 & \mbox{otherwise}, 
\end{array}
\right.
\eea
where $\hz:=-q^{-1}bc$. 
One can think of this cocycle  as the ``Chern cyclic cocycle,"
and the invariants it computes as the Chern numbers of quantum vector bundles.
The fact that it is in degree zero is a quantum effect
caused by the non-classical structure of $HC^*(\asq)$ (see~\cite{mnw91}). 
In the classical case the corresponding cocycle is in degree two,
 as it comes from the volume form of the two-sphere.

The pairing $\<[\tau^1],[e_{-1}]\>=-1$, $\<[\tau^1],[1]\>=0$ 
can be used  to conclude the non-cleftness of the quantum principal
Hopf fibration~\cite[Corollary~4.2]{hm98}. Furthermore, taking advantage of
the linearity of the pairing and  the above equalities, we can also conclude
that $\IZ\oplus\IZ\inc K_0(\asq)$. On the other hand, one can directly check that
the other cyclic cocycle~$\tau^0$ pairs unitally with 
above projectors:
\beq\Label{rank}
\langle [\tau^0],[e_{\mu}]\rangle=(\he\ci Tr)(e_\mu)=1,~~~\mu\in\IZ.
\eeq
We can think of $\tau^0$,
$\tau^1$ as a (possibly incomplete) coordinate system for~$K_0(\asq)$ determining
the rank and Chern number respectively. The point of this paper is that 
for any $\mu\in\IZ$ there exists
a rank one projector matrix (quantum line bundle) with its Chern number equal 
to~$\mu$.

\section{Chern-Connes pairing for quantum Hopf line bundles}

We are to compute the pairing between the Chern cyclic cocycle
$\tau^1$ and both left and right projector matrices of quantum Hopf 
line bundles~$P_\mu$.
We refer to the thus obtained invariants as the left and right Chern numbers 
respectively. 
Since $\tau^1$ is a 0-cyclic cocycle, the pairing is given by the formula
$\<[\tau^1],[p]\>=(\tau^1\ci Tr)(p)$, where $p\in M_n(\asq)$, $p^2=p$, and
$Tr: M_n(\asq)\ra\asq$ is the usual matrix trace. We have:
\bth\Label{k}
The left Chern number and the winding number  of any quantum Hopf line bundle
coincide: $(\tau^1\ci Tr)(e_\mu)=\mu,\;\mu\in\IZ$.
\ethe\bpf
We need to consider two cases: $\mu < 0$ and $\mu > 0$. (The case $\mu=0$ is
evident, as $e_0=1$ and $\tau^1$ annihilates numbers.)

{\em Case $\mu < 0$:~} To simplify notation put $n=-\mu$. Let us first compute
the trace of $e_{-n}$ (see~(\ref{ekl})) as a polynomial in~$\hz:=-q^{-1}bc$:
\begin{eqnarray*}
Tr(e_{-n})\!\!\!\!\!\! &&
=\sum_{k=0}^{n}{\scriptsize\pmatrix{n\cr k\cr}}_{\!q^2}
(-q)^{k}a^{n-k}c^kb^{k}d^{n-k}
\\ &&
=\sum_{k=0}^{n}{\scriptsize\pmatrix{n\cr k\cr}}_{\!q^2}
(-q)^{k}q^{-2k(n-k)}(bc)^{k}a^{n-k}d^{n-k}
\\ &&
=\sum_{k=0}^{n}{\scriptsize\pmatrix{n\cr k\cr}}_{\!q^2}
q^{-2k(n-k-1)}\hz^{k}\prod^{n-k-1}_{l=0}(1-q^{-2l}\hz).
\end{eqnarray*}
The last step follows from the quantum determinant formula $ad=1+q^{-1}bc=1-\hz$
(see~(\ref{comm})) and standard induction. 
(The expression $\prod^{-1}_{l=0}(\cdots)$ is understood
as~1.) To apply $\tau^1$ to $Tr(e_{-n})$ we need to know more explicitly the
coefficients of the above polynomial. For this reason let us recall the definition
of {\em shifted binomials} (cf.~\cite[p.173]{mnw91}). Let $k[t]$ denote a polynomial
ring in one variable and $x\in k$. For natural numbers $0\leq l\leq\nu$ we define
the $x$-shifted binomial \mbox{\scriptsize$\left[\matrix{\nu\cr l\cr}\right]_{\!x}$} 
by the equality
\[
\sum^{\nu}_{l=0}\mbox{\scriptsize$\left[\matrix{\nu\cr l\cr}\right]_{\!x}$}t^l
:=\prod^{\nu-1}_{l=0}(1+x^lt)
\]
Now we can use the above calculations and (\ref{tr}) to compute 
the Chern-Connes pairing between $[\tau^1]$ and $[e_{-n}]$:
\begin{eqnarray*}
\langle [\tau^1],[e_{-n}]\rangle\!\!\!\!\!\! &&
=(\tau^1\ci Tr)(e_{-n})
\\ &&
=\tau^1\left(\sum_{k=0}^{n}{\scriptsize\pmatrix{n\cr k\cr}}_{\!q^2}
q^{-2k(n-k-1)}\hz^{k}\sum^{n-k}_{l=0}
\mbox{\scriptsize$\left[\matrix{n-k\cr l\cr}\right]_{\!q^{-2}}$}(-\hz)^l\right)
\\ &&
=\tau^1\left(\sum_{k=0}^{n}{\scriptsize\pmatrix{n\cr k\cr}}_{\!q^2}
q^{-2k(n-k-1)}\hz^{k}\sum^{n}_{m=k}
\mbox{\scriptsize$\left[\matrix{n-k\cr m-k\cr}\right]_{\!q^{-2}}$}(-\hz)^{m-k}\right)
\\ &&
=\tau^1\left(\sum_{m=0}^{n}(-\hz)^{m}
\sum^{m}_{k=0}{\scriptsize\pmatrix{n\cr k\cr}}_{\!q^2}
q^{-2k(n-k-1)}(-1)^{k}
\mbox{\scriptsize$\left[\matrix{n-k\cr m-k\cr}\right]_{\!q^{-2}}$}\right)
\\ &&
=\sum_{m=1}^{n}(1-q^{2m})^{-1}(-1)^{m}
\sum^{m}_{k=0}{\scriptsize\pmatrix{n\cr k\cr}}_{\!q^2}
q^{-2k(n-k-1)}(-1)^{k}
\mbox{\scriptsize$\left[\matrix{n-k\cr m-k\cr}\right]_{\!q^{-2}}$}\, .
\end{eqnarray*}
The point is to prove that the just computed number equals~$-n$. To do so let
us assume for the time being that $q\in (0,1)$, so that we can use the $C^*$-algebraic
framework. Recall that the $0$-cyclic cocycle $\tau^1$ comes from a 1-summable
Fredholm module over~\asq~\cite[p.175]{mnw91}. Hence, by the Noncommutative Index
Theorem, the Chern-Connes pairing between $[\tau^1]$ and any element of $K_0(\asq)$
is necessarily an integer --- the index of an appropriate Fredholm operator. 
(See, e.g., \cite[p.297]{c-a94},\cite[p.54]{c-a88},\cite[Section~12.2.5]{l-jl97}.)
Thus we have 
\[
(\tau^1\ci Tr)(e_{-n})=
\sum_{m=1}^{n}(1-q^{2m})^{-1}(-1)^{m}
\sum^{m}_{k=0}{\scriptsize\pmatrix{n\cr k\cr}}_{\!q^2}
q^{-2k(n-k-1)}(-1)^{k}
\mbox{\scriptsize$\left[\matrix{n-k\cr m-k\cr}\right]_{\!q^{-2}}$}
\in\IZ
\]
for any~$q\in (0,1)$. Observe now that, since $(\tau^1\ci Tr)(e_{-n})$ is a rational
function of $q$, it is continuous whence constant on $(0,1)$, by the connectedness
of $(0,1)$ and the above integrality property. The only rational function on
$\IC\setminus\{0, \mbox{roots of 1}\}$ which is constant on the open interval $(0,1)$
is a constant function. Therefore $(\tau^1\ci Tr)(e_{-n})$ is independent of 
$q\in(1,\infty)$,  and we have 
$(\tau^1\ci Tr)(e_{-n})=\lim_{q\ra\infty}(\tau^1\ci Tr)(e_{-n})$. Now it suffices
to show that 
\[
\lim_{q\ra\infty}~(1-q^{2m})^{-1}(-1)^{m}
\sum^{m}_{k=0}{\scriptsize\pmatrix{n\cr k\cr}}_{\!q^2}
q^{-2k(n-k-1)}(-1)^{k}
\mbox{\scriptsize$\left[\matrix{n-k\cr m-k\cr}\right]_{\!q^{-2}}$}=-1,~~~
m,\, n>0.
\]
To this end we need to analyse the asymptotic behaviour of the fractions
\begin{eqnarray*}
&& F_{n,m,k}(q):={\scriptsize\pmatrix{n\cr k\cr}}_{\!q^2}
q^{-2k(n-k-1)}(-1)^{k}
\mbox{\scriptsize$\left[\matrix{n-k\cr m-k\cr}\right]_{\!q^{-2}}$}
\\ && =
\mbox{$
\frac{(-1)^{k}(q^2-1)...(q^{2n}-1)}
{q^{2k(n-k-1)}(q^2-1)...(q^{2k}-1)(q^2-1)...(q^{2(n-k)}-1)}
\left(\frac{1}{q^{(m-k-1)(m-k)}}+...+\frac{1}{q^{(n-k+1)(n-k)-(n-m+1)(n-m)}}\right)
$},
\end{eqnarray*}
where the sum in parenthesis 
equals~\mbox{\scriptsize$\left[\matrix{n-k\cr m-k\cr}\right]_{\!q^{-2}}$}. 
The dominating term of $F_{n,m,k}(q)$ if proportional to $q$ to the power of
$n(n+1)-k(k+1)-(n-k)(n-k+1)-2k(n-k-1)-(m-k-1)(m-k)=-k^2+k(1+2m)+m-m^2$. Thinking
of this expression as a function of~$k$, we see that it is biggest for $k=m$.
Hence the dominating term in $\sum_{k=0}^{m}F_{n,m,k}(q)$ is proportional to~$q^{2m}$.
All other terms are proportional to $q$ to some power strictly less than $2m$,
and will be vanished by $(1-q^{2m})^{-1}$. Therefore
\begin{eqnarray*}
&& \lim_{q\ra\infty}~(1-q^{2m})^{-1}(-1)^{m}\sum^{m}_{k=0}F_{n,m,k}(q)
\\ && =
\lim_{q\ra\infty}~(1-q^{2m})^{-1}(-1)^{m}
\frac{(-1)^{m}(q^2-1)...(q^{2n}-1)}
{q^{2m(n-m-1)}(q^2-1)...(q^{2m}-1)(q^2-1)...(q^{2(n-m)}-1)}
\\ && =
-\lim_{q\ra\infty}
\frac{(q^2-1)...(q^{2n}-1)}
{(q^{2m}-1)q^{2m(n-m-1)}(q^2-1)...(q^{2m}-1)(q^2-1)...(q^{2(n-m)}-1)}
\\ && = -1.
\end{eqnarray*}
This proves that $(\tau^1\ci Tr)(e_{-n})=\sum_{m=1}^{n}(-1)=-n$, as needed.

{\em Case $\mu>0$:}~ The reasoning is similar to that of the previous case, though
the calculation of the limit is more straightforward. Put $n=\mu$. First we compute:
\begin{eqnarray*}
Tr(e_{n})\!\!\!\!\!\! &&
=\sum_{k=0}^{n}{\scriptsize\pmatrix{n\cr k\cr}}_{\!q^2}
(-q)^{-k}b^{k}d^{n-k}a^{n-k}c^k
\\ &&
=\sum_{k=0}^{n}{\scriptsize\pmatrix{n\cr k\cr}}_{\!q^2}
\hz^{k}\prod^{n-k-1}_{l=0}(1-q^{-2(l+1)}\hz)
\\ &&
=\sum_{k=0}^{n}{\scriptsize\pmatrix{n\cr k\cr}}_{\!q^2}
\hz^{k}\sum^{n-k}_{l=0}
\mbox{\scriptsize$\left[\matrix{n-k\cr l\cr}\right]_{\!q^{2}}$}(-q^{2}\hz)^l
\\ &&
=\sum_{k=0}^{n}{\scriptsize\pmatrix{n\cr k\cr}}_{\!q^2}\hz^{k}
\sum^{n}_{m=k}
\mbox{\scriptsize$\left[\matrix{n-k\cr m-k\cr}\right]_{\!q^{2}}$}(-q^{2}\hz)^{m-k}
\\ &&
=\sum_{m=0}^{n}\hz^{m}
\sum^{m}_{k=0}{\scriptsize\pmatrix{n\cr k\cr}}_{\!q^2}
(-1)^{m-k}q^{2(m-k)}
\mbox{\scriptsize$\left[\matrix{n-k\cr m-k\cr}\right]_{\!q^{2}}$}.
\end{eqnarray*}
Using the Noncommutative Index Theorem the same way as before, we can again conclude
that $(\tau^1\ci Tr)(e_{n})\in\IZ$ for $q\in(0,1)$. Also by the same argument as 
before, the integrality of $(\tau^1\ci Tr)(e_{n})$ entails that it is independent of
$q\in(0,1)$. We can therefore compute it
by taking the limit: 
\begin{eqnarray*}
(\tau^1\ci Tr)(e_{n})\!\!\!\!\!\! 
&& =
\lim_{q\ra 0}~(\tau^1\ci Tr)(e_{n})
\\ && =
\lim_{q\ra 0}~\sum_{m=1}^n
(1-q^{2m})^{-1}\sum_{k=0}^m{\scriptsize\pmatrix{n\cr k\cr}}_{\!q^2}
(-1)^{m-k}q^{2(m-k)}
\mbox{\scriptsize$\left[\matrix{n-k\cr m-k\cr}\right]_{\!q^{2}}$}
\\ && =
\sum_{m=1}^n\sum_{k=0}^m(-1)^{m-k}\hd_{km} = n.
\end{eqnarray*}
\epf 
\bre\Label{2}\em
It follows from a direct calculation that 
$(\tau^1\ci Tr)(e_{-2})=(-1-q^{-2})+(-1+q^{-2})=-2$,
where the first and the second term correspond to $m=1$ and $m=2$ respectively
in the general sum. Similarly, $(\tau^1\ci Tr)(e_{2})=(1+q^{2})+(1-q^{2})=2$. 
\ere
As for the right projective structure of $P_\mu$, 
one can infer directly from formulas (\ref{ekl}), 
(\ref{flk}) and Theorem~\ref{k} that
\[
(\tau^1\ci Tr)(f_\mu)=(\tau^1\ci Tr)(e_{-\mu})=-\mu.
\]
Hence we have:
\bco\Label{bi}
The right Chern number of any quantum Hopf line bundle equals the opposite of its
winding number, i.e., $(\tau^1\ci Tr)(f_\mu)=-\mu,\;\mu\in\IZ$.
\eco
Let us now consider further consequences of Theorem~\ref{k}. Note first that
due to the direct sum decomposition $\slq=\bigoplus_{\mu\in\IZ}P_\mu$, one can
say that the coordinate ring of $SL_q(2)$ decomposes into {\em mutually 
$K_0$-non-equivalent} left and right projective finitely generated \asq-bimodules.
As for the structure of $K_0(\asq)$, Theorem~\ref{k} provides us with
an estimate of its positive cone.  Indeed, combining (\ref{rank}) with 
 Theorem~\ref{k} yields:
\bco\Label{k+}
The image of the positive cone of $K_0(\asq)$ under 
$(\tau^0,\tau^1):K_0(\asq)\lra\IZ\times\IZ$
contains $\IN\times\IZ$.
\eco


{\bf Acknowledgements:} The author was partially supported 
by  the  NATO and CNR postdoctoral fellowships and KBN grant
2 P03A 030 14. 


\end{document}